\documentclass[a4paper]{article}

\usepackage[english]{babel}

\usepackage{ifxetex}
\ifxetex
\usepackage{fontspec}
\else
\usepackage[utf8]{inputenc}
\fi

\usepackage[a4paper,top=3cm,bottom=2cm,left=3cm,right=3cm,marginparwidth=1.75cm]{geometry}

\usepackage{amsmath, amssymb, amsthm}
\usepackage{graphicx}
\usepackage[colorinlistoftodos]{todonotes}
\usepackage{hyperref}
\usepackage{bbold}

\usepackage[style=numeric-comp,backend=bibtex]{biblatex}
\newtheorem{theorem}{Theorem}

\newcommand\kk{\mathbb{k}}
\newcommand\id{\mathbb 1}
\newcommand\img{\operatorname{img}}

\newcommand\Hom{\operatorname{Hom}}

\newcommand\Ext{\operatorname{Ext}}
\newcommand\Dgm{\operatorname{Dgm}}

\newcommand\into{\hookrightarrow}
\newcommand\isoto{\xrightarrow{\sim}}

\newcommand\XX{\mathbb{X}}
\newcommand\VV{\mathbb{V}}
\newcommand\WW{\mathbb{W}}

\title{Algorithms in $A_\infty$-algebras}
\author{Mikael Vejdemo-Johansson}

\begin{document}
\maketitle
\textit{To Tornike Kadeishvili}
\begin{abstract}
  Building on Kadeishvili's original theorem inducing $A_\infty$-algebra structures on the homology of dg-algebras, several directions of algorithmic research in $A_\infty$-algebras have been pursued.
  In this paper we will survey work done on calculating explicit $A_\infty$-algebra structures from homotopy retractions; in group cohomology; and in persistent homology.
\end{abstract}

\section{Introduction}

In his 1980 paper~\cite{kadeishvili80}, Kadeishvili proved that the homology of any dg-algebra has an induced $A_\infty$-algebra structure.
The proof itself is by induction on the arity of the operation, writing out explicitly how to create $m_n(a_1,\dots,a_n)$ in terms of lower order operations.
By giving these expressions, the proof is basically a formulation of what an algorithm for the computation of an $A_\infty$ algebra structure on the homology of a dg-algebra would look like.

Even though this initial setup has a strongly algorithmic flavor, actual algorithms and computer-supported calculations of $A_\infty$-algebra structures emerged far later. 
In this paper, we plan to give an overview of work on creating software computing $A_\infty$-algebra structures and the contexts and techniques used for these.
Merkulov~\cite{merkulov_sha_kahler99} has given Kadeishvili's construction a more concrete, and combinatorially accessible presentation -- but the focus for this paper is in explicit algorithmic computation with computer implementations, preferably publically accessible.

\section{Background}

We fix a ring $\kk$. All tensor products are over $\kk$ unless otherwise noted, and tensor powers are denoted by $A^{\otimes n} = A\otimes\dots\otimes A$.

A graded $k$-vector space $A$ is an \emph{$A_\infty$-algebra} if one of the
following equivalent conditions hold

\begin{enumerate}
\item There is a family of maps $\mu_i\colon A^{\otimes i}\to A$, called higher multiplications fulfilling the Stasheff identities 
  \[
  \operatorname{St}_n\colon \sum_i\sum_j \mu_i\circ_j\mu_{n-i} = 0\quad.
  \]
\item There is a family of chain maps from the cellular chain complex of the associahedra to appropriate higher endomorphisms of $A$
  \[
  \mu_n\colon C_*(K_n)\to\operatorname{Hom}(A^{\otimes n}, A)\quad.
  \]
\item $A$ is a representation of the free dg-operad resolution $\mathcal Ass_\infty$ of the associative operad. 
\end{enumerate}

An \emph{$A_\infty$-coalgebra} is defined by dualizing this definition.

The structure was introduced by Jim Stasheff in \cite{stasheff_ha_hspaces63}.
Good introductory surveys have been written by Lu, Palmieri, Wu and Zhang \cite{lu_wu_palmieri_zhang_Aoo_ring04,lu_wu_palmieri_zhang_Aoo_Ext06} as well as by Bernhard Keller \cite{keller_intro01,keller_Aoo_repr00}.

A graded $k$-vector space $A$ is a differential graded algebra (dg-algebra) if it is equipped with a differential operator $\partial: A\to A$ of degree $-1$ and an associative multiplication $m_2: A\otimes A\to A$ of degree $0$, such that the Leibniz rule holds:
\[
\partial m_2(x, y) = m_2(\partial x, y) + (-1)^{|x|}m_2(x, \partial y)
\]
A module over a dg-algebra $A$ is a graded vectorspace $M$ with a differential operator $\partial: M\to M$ and an associative multiplication $m_2: M\otimes A\to M$ that obeys the Leibniz rule.

Kadeishvili proved in his 1980 paper~\cite{kadeishvili80} that the homology of a dg-algebra has an inherited and quasi-isomorphic $A_\infty$-algebra structure. The proof starts out defining $\mu_1=0$ and $\mu_2([x_1]\otimes [x_2]) = (-1)^{|x_1|+1}[x_1\cdot[x_2]$, as well as starting to define an $A_\infty$-morphism $f$ by $f_1$ simply a cycle-choosing homomorphism.

The proof proceeds by induction\footnote{which translates well to recursion: this is how algorithms enter the picture}: if all $m_j$ and $f_j$ have been defined for $j<i$, then let
\[
U_n = \sum_{s=1}^n m_2 \circ f_s \otimes f_{n-s}) + 
\sum_{k=0}^{n-2}\sum_{j=2}^{n-1} f_{n-j+1} \circ \id^{\otimes k}\otimes m_j \otimes \id^{\otimes n-j-k+1}
\]

The Stasheff axioms can be translated to 
\[
m_1\circ f_n = (f_1\circ m_n - U_n)
\]

The right hand side is homological to zero, so we can pick $f_n$ to be a bounding element of this difference.
Extend by linearity.

The structure of this argument lends itself excellently well to concrete and algorithmic calculations, and there has been a few approaches to algorithmic and computer-aided $A_\infty$-algebra work. 
There are three main topics that emerge, and we will dedicated a chapter to each of them. 

First up, in Section~\ref{sec:reductions}, we will review Ainhoa Berciano's work on contractions of dg-algebras to dg-modules, with implementations in the computer algebra system Kenzo.

Next, in Section~\ref{sec:group-cohomology} we will go to the realm of group cohomology. 
Mikael Vejdemo-Johansson worked on algorithms to directly calculate $A_\infty$-algebra structures on the modular group cohomology of $p$-groups, and generated a number of ways to recognize feasibility of the calculation as well as a stopping criterion. 
Stephan Schmid, later on, used some of Vejdemo-Johansson's results in a concrete calculation of the $A_\infty$-algebra structure on the modular group cohomology of the symmetric group $S_p$. 

Finally, in Section~\ref{sec:persistent-infinity}, we will describe recent work by Murillo and Belchí on using $A_\infty$-coalgebra structures on persistent homology rings to create new perspectives on bottleneck distances and stability of persistence barcodes.

\section{Reductions}
\label{sec:reductions}

Ainhoa Berciano's work~\cite{berciano-2006, berciano2010computational, berciano2012searching, berciano2004coalgebra, berciano2010coalgebra} starts with perturbation theory.
This framework has as its core result the Basic Perturbation Lemma~\cite{brown1967twisted}, that describes how a \emph{contraction} changes under perturbation.

A \emph{contraction} connects two dg-modules $M$ and $N$, abstracting the homotopy concepts of deformation retracts. 
A contraction consists of morphisms $f: N\to M$, $g: M\to N$ and $\phi: N\to N$ such that $f$ and $g$ are almost an isomorphism -- up to a homotopy operation in $N$. In other words, we require
\[
fg = \id_M
\quad
gf + \phi\partial_N + \partial_N\phi = \id_N
\quad
f\phi = 0
\quad
\phi g = 0
\quad
\phi\phi = 0
\]

A contraction preserves homology: $H(N)$ is canonically isomorphic to $H(M)$, but this isomorphism tends not to transfer algebraic structures from $N$ to $M$.

If the differential structure on $N$ is perturbed: instead of boundary operator $\partial_N$, $N$ is equipped with a new boundary operator $\partial_N+\delta$, then the Basic Perturbation Lemma produces a new contraction.
The requirement for this construction is that $\phi\delta$ is pointwise nilpotent: for any $x\in N$ there is an $n$ so that $(\phi\delta)^n(x)=0$.

Then there is a new contraction $f_\delta, g_\delta, \phi_\delta$ between $N$ equipped with the boundary operator $\partial_N+\delta$ and $M$ equipped with the boundary operator $\partial_M+\partial_\delta$ given by
\begin{align*}
\partial_\delta &= f\delta\sum_{i\geq 0}(-1)^i(\phi\delta)^i g  &
f_\delta &= f\left(1-\delta\sum_{i\geq 0}(-1)^i(\phi\delta)^i\phi\right)
\\
g_\delta&=\sum_{i\geq 0}(-1)^i(\phi\delta)^ig &
\phi_\delta&=\sum_{i\geq 0}(-1)^i(\phi\delta)^i\phi
\end{align*}

Berciano generates an algorithm for transferring $A_\infty$-coalgebra structures between DG-modules using the \emph{tensor trick}~\cite{gugenheim-lambe-stasheff-91}.
The tensor trick starts with a dg-coalgebra $C$, a dg-module $M$ and a contraction from $C$ to $M$.
Take the tensor module of the desuspension of all components in this contraction to produce a new contraction.
With the cosimplicial differential, we can use the Basic Perturbation Lemma and obtain a new contraction.
The tilde cobar differential generates an induced $A_\infty$-coalgebra structure, explicitly given by the comultiplication operations
\begin{equation}\label{eqn:Aoo-coalgebra}
\Delta_i=(-1)^{[i/2]+i+1}f^{\otimes i}\Delta^{[i]}\phi^{[\otimes(i-1)]}\dots\phi^{[\otimes 2]}\Delta^{[2]}g ;
\qquad
\Delta^{[k]} = \sum_{i=0}^{k-2}(-1)^i\id^{\otimes i}\otimes\Delta\otimes\id^{k-i-2}
\end{equation}

This derivation allows Berciano to prove~\cite{berciano2004coalgebra} that in $H_*(K(\pi, n); \mathbb{Z}_p)$ for a finitely generated abelian group $\pi$, the only non-null morphisms in the $A_\infty$-coalgebra structure have to have order $i(p-2)+2$ for some non-negative integer $i$. 

The final formula in~\ref{eqn:Aoo-coalgebra} is concrete enough that it has been implemented on the computer algebra platform Kenzo in the packages ARAIA (\textbf{A}lgebra \textbf{R}eduction \textbf{A}-\textbf{I}nfinity \textbf{A}lgebra) and CRAIC (\textbf{C}oalgebra \textbf{R}eduction \textbf{A}-\textbf{I}nfinity \textbf{C}oalgebra).

\section{Group Cohomology}
\label{sec:group-cohomology}

Fix a group $G$ and a field $\kk$. The cohomology algebra of the Eilenberg-MacLane space $H^*(K(G, 1))$ is isomorphic to the Ext-algebra $\Ext_{\kk G}(\kk, \kk)$ of the group ring and is called the group cohomology $H^*(G)$.
Because of the connection to the Ext-algebra, the group cohomology can be calculated from the composition dg-algebra of $\Hom(F_*, F_*)$ for a free resolution $F_*\to\kk$ in the category of $G$-modules.
Several computer algebra systems, including Magma~\cite{magma} and GAP~\cite{GAP4} support calculations with $G$-modules.
In such a system, we create a free resolution $F_*$ of $\kk$.
Chain maps $F_*\to F_*$ are then represented by a sequence of maps, one for each degree, each determined by lower-dimensional maps through commutativity of the corresponding squares in the chain map diagram.
With $\Hom(F_*, F_*)$ represented, we can compute $H^*G$ as $H_*\Hom(F_*, F_*)$. 

Since the $\Hom(F_*, F_*)$ is a dg-algebra, by Kadeishvili's theorem, $H^*G$ has an induced $A_\infty$-algebra structure.

Vejdemo-Johansson~\cite{mikAooCncm07,vejdemo2010blackbox,vejdemo2008computation} studies this $A_\infty$-algebra structure from a strictly algorithmic perspective.

\subsection{Blackbox computation of $A_\infty$}
\label{sec:blackb-comp-a_infty}

A cornerstone of Vejdemo-Johansson's approach to computing $A_\infty$-algebras is the following theorem (\cite[Theorem 3]{vejdemo2010blackbox}):

\begin{theorem}\label{thm:mvj:linearity}
  If $A$ is a dg-algebra and 
  \begin{enumerate}
  \item There is an element $z\in H_*A$ generating a polynomial subalgebra (ie is not a torsion element)
  \item $H_*A$ is a free $\kk[z]$-module
  \item\label{thm:cond:commutativity} $H_*A$ has a $\kk[z]$ linear $A_{n-1}$-algebra structure induced by the dg-algebra structure on $A$, such that $f_1(z)f_k(a_1, \dots, a_k) = f_k(a_1, \dots, a_k)f_1(z)$
  \item We have a chosen $\kk[z]$-basis $b_1, \dots$ of $H_*A$ and all $m_k(v_1, \dots, v_k)$ and $f_k(v_1, \dots, v_k)$ are chosen by Kadeishvili's algorithm for all combinations of basis elements $v_j\in\{b_1,\dots\}$
  \end{enumerate}

  Then a choice of $m_n$ and $f_n$ by Kadeishvili's algorithm for all input values taken from this $\kk[z]$-basis extends to a $\kk[z]$-linear $A_n$-algebra structure on $H_*A$ induced by the dg-algebra structure on $A$.
\end{theorem}

The condition~\ref{thm:cond:commutativity} says that for the $A_\infty$-morphism $H_*A\to A$ produced by Kadeishvili's construction, the cycle chosen for $z$ commutes -- \textbf{on a chain level} -- with each chain map chosen for the higher operations.
This is the key condition for the theorem -- and also the one that makes the theorem most fragile.

The theorem tells us we can construct an $A_\infty$-algebra structure step by step.
If there is one of these non-torsion central elements $z$, we can reduce the complexity of $H_*A$ for the purpose of calculating its higher operations -- if we find a family of central elements $z_1, \dots, z_k$ such that $H_*A$ is a finite module over $\kk[z_1,\dots,z_k]$, then it is enough to study the finitely many basis elements $b_1,\dots,b_m$ in a presentation of $H_*A$ as a  $\kk[z_1,\dots,z_k]$-module.
This makes each calculation a finite (though large) in terms of the number of input combinations that need to be studied.
Using this theorem is easier if -- as is the case for cyclic groups -- the resolution $F_*$ is periodic.

Theorem~\ref{thm:mvj:linearity} makes it easier to extend from $A_{n-1}$ to $A_n$, using a  condition that can be checked for each extension step.
Once the condition -- commutativity of the representative chain maps -- fails, the structure calculated thus far is valid, but further extensions are obstructed.
The key to bring computational effort down to a finite time endeavour lies in \cite[Theorem 5]{vejdemo2010blackbox}:

  \begin{theorem}
    Let $A$ be a dg-algebra.
    If in an $A_{2q-2}$-algebra structure on $H_*A$, $f_k=0$ and $m_k=0$ for all $q\leq k\leq 2q-2$ then the $A_{2q-2}$-structure is already an $A_\infty$-structure with all higher $f_n$ and all higher $m_n$ given by zero maps.
  \end{theorem}

Through finding central elements, the infinitely many basis elements of $H_*A$ can be brought down to a finite number of basis elements to check.
And by finding a large enough gap, in which all chain representations and all products vanish, the computation can be terminated producing a result.

This approach was implemented as a module distributed with Magma~\cite{magma}, and was used both to confirm Madsen's~\cite{madsen_phd} computation of $A_\infty$-algebra structures on the group cohomology of cyclic groups and to conjecture~\cite{vejdemo2008computation} the start of an $A_\infty$-structure on the cohomology on some dihedral groups.

\subsection{The Saneblidze-Umble diagonal}
\label{sec:sanebl-umble-diag}

In~\cite{saneblidze-umble-2004}, Saneblidze and Umble gave an explicit construction for a diagonal on the associahedra.
This construction translates directly to a method to combine $A_\infty$-algebra structures on $V$ and $W$ into an $A_\infty$-algebra structure on $V\otimes W$.

Vejdemo-Johansson uses this construction in~\cite{mikAooCncm07} to prove non-triviality of some operations on $H^*(C_n\times C_m)$.
From results by Berciano and Umble~\cite{berciano2011some}, we know that any non-trivial operation on this group cohomology of arity less than $n+m-1$ has to have arity $2, n, m$ or $n+m-2$.
Berciano also shows~\cite{berciano-2006} that any non-zero higher coproduct on $H_*(C_q\times C_q)$ has arity $k(q-2)+2$ for some $k$. 
In addition to the induced operations in arities $2, n, m$ and $n+m-2$, Vejdemo-Johansson shows that there are non-trivial operations of arity $2n+m-4$ and $n+2m-4$.
The original article states a far more generous claim: that all the arities $k(n-2) + k(m-2)+2$, $(k-1)(n-2) + k(m-2) + 2$ and $k(n-2) + (k-1)(m-2) + 2$ have non-zero operations -- this argument turned out to have a subtle flaw, and was retracted.
More details are available in~\cite{vejdemo2008computation}.

Any practical use of the Saneblidze-Umble diagonal would benefit greatly from a computer-facilitated access to the coefficients of the diagonal construction.
In an unpublished preprint~\cite{vejdemo2007enumerating}, Vejdemo-Johansson provides a computer implementation of an algorithm to enumerate the Saneblidze-Umble terms.

\subsection{Symmetric groups}
\label{sec:symmetric-groups}

Schmid~\cite{Schmid2014} studies the group cohomology of the symmetric group $S_p$ on $p$ elements, with coefficients in the finite field $\mathbb{F}_p$ with $p$ elements.
For this group cohomology, he presents a basis with which he is able to prove that the only non-trivial $A_\infty$-operations on $H^*S_p$ are of arity $2$ and $p$.
To do this, he goes through large and somewhat onerous explicit calculations to show that there is a periodic projective resolution of $\mathbb{F}_p$ over $\mathbb{F_p}S_p$, and that the resolution has a large enough gap to allow the use of Vejdemo-Johansson's theorem.

\section{Persistent A-infinity}
\label{sec:persistent-infinity}

Persistent homology and cohomology form the cornerstone of the fast growing field of Topological Data Analysis.
The fundamental idea is to study the homology functor applied to diagrams of topological spaces
\[
\VV : \VV_0 \into \VV_1 \into \dots
\]

These spaces are often generated directly from datasets, by constructions such as the \v{C}ech construction: for data points $\XX = \{x_0, \dots, x_N\}$, an abstract simplical complex $\check{C}_\epsilon$ has as its vertices $\XX$ and includes a simplex $[x_{i_0}, \dots, x_{i_d}]$ precisely if the intersection of balls $\bigcap_{j=0}^d B_\epsilon(x_{i_j})$ is non-empty.
If $\epsilon$ increases, no intersections will become empty, and so no simplices will vanish.
So the \v{C}ech complexes, as $\epsilon$ sweeps from 0 to $\infty$, generates a nested sequence of topological spaces.

The inclusion maps $\iota_i^j:\VV_i\to\VV_j$ lift by functoriality to linear maps on homology: 
$H_*(\iota_i^j): H_*\VV_i\to H_*\VV_j$.
We may define a persistent homology group as the image $PH_*^{i,j}(\VV) = \img H(\iota_i^j)$.

For more details on the data analysis side, we recommend the surveys \cite{carlsson2009topology,ghrist2008barcodes,vejdemo2014sketches}

\subsection{Barcodes and stability}
\label{sec:barcodes-stability}

As the homology functor is applied to the diagram of topological spaces, using coefficients from a field $\kk$ for the homology computation, the result is a diagram of vector spaces.
By either imbuing the resulting diagram with the structure of a module over the polynomial ring $\kk[t]$, or as representations of a quiver $Q$ of type $A_n$, the corresponding classification theorems produce a decomposition of $H_*(\VV)$ into a direct sum of \emph{interval modules}.
These interval modules are $\mathbb N$-graded modules defined by a pair of indices $b, d$, and are defined as $0$-dimensional for degrees $k<b$ and for degrees $k > d$.
For degrees from $b$ to $d$, the interval module is one-dimensional, with identity maps connecting each space to the next.

Thus, the homology of a diagram $\VV$ of topological spaces with field coefficients can be described by a multiset $\Dgm(\VV) = \{ (b_i, d_i) \}_{i\in I}$, called the \emph{persistence barcode} or \emph{persistence diagram} of the diagram $\VV$. 
The dimension of the persistent homology group $PH_*^{i,j}(\VV)$ is exactly the number of intervals $(b_k, d_k)$ in $\Dgm(\VV)$ such that $b_k \leq i \leq j \leq d_k$.

Between any pair of such diagrams we can create a distance called the \emph{bottleneck distance}. Setting $\Delta = \{(x,x) : x\in\mathbb R\}$, this distance is defined by:

\[
d_B(\Dgm(\VV, \WW)) = \inf_{\gamma: \VV\cup\Delta\isoto\WW\cup\Delta} \quad\max_{v\in\VV\cup\Delta} d(v, \gamma v)
\]

This distance measures the largest displacement needed to change $\Dgm(\VV)$ into $\Dgm(\WW)$, while allowing intervals to disappear into and emerge from the infinite set of possible 0-length intervals.

First with the bottleneck distance (and in later research more sophisticated), a range of \emph{stability theorems} have been proven, starting in \cite{cohen2007stability}. Good overviews can be found in \cite{vejdemo2014sketches,chazal2012structure}.
These theorems take the shape of 

\begin{theorem}[Stability meta-theorem]
  If $d(\VV, \WW) < \epsilon$, then $d'(\Dgm(H(\VV)), \Dgm(H(\WW))) < \epsilon$ for specific choices of distances $d$ and $d'$.
\end{theorem}

\subsection{$A_\infty$ in persistence}
\label{sec:a_infty-persistence}

Murillo and Belchí introduced in \cite{belchi2015ainfty,belchi2017optimising} an $A_\infty$-coalgebra approach to barcode distances.
If each new cell introduced in the step from $\VV_i$ to $\VV_{i+1}$, with all the $\VV_j$ chosen to be CW-complexes, and working over the rationals $\mathbb{Q}$, then there is a set of compatible choices of $A_\infty$-coalgebras for the entire sequence.

They define a $\Delta_n$-persistence group 
\[
\Delta_nPH_*^{i,j}(\VV) = \img (H_*(\iota_i^j)|_{\bigcap_{k=i}^j\ker(\Delta_n^k\circ\iota_i^k}
\]

In other words, the $\Delta_n$-persistence group retains from the ordinary persistence groups precisely those elements out of $H_*(\VV_i)$ whose images in each $\VV_k$ vanish under application of the higher coproduct $\Delta_n$.

These $\Delta_n$-persistence groups generate $\Delta_n$-persistence barcodes as multisets $\Dgm_{\Delta_n}(\VV)$ of intervals $[b,d]$. From these barcodes, the dimension of $\Delta_nPH_*^{i,j}(\VV)$ equals the number of $[b_k, d_k]\in\Dgm_{\Delta_n}(\VV)$ such that $b_k\leq i\leq j\leq d_k$.
As pointed out by the authors, these higher order barcodes may ``flicker'' in a way that classical persistence strictly avoids: the same element can exist over several disjoint intervals.
This has been carefully avoided in the greater literature on persistent homology: the flickering behavior invites wild representation theories, where the decomposition that generates barcodes is no longer available.

\subsection{$A_\infty$ bottleneck distance}
\label{sec:a_infty-bottl-dist}

Herscovich \cite{herscovich2014higher} introduces a novel metric on persistent homology.
Herscovich constructs a metric on locally finite Adams graded minimal $A_\infty$-algebras, and then quotients by quasi-isomorphism to establish a metric on persistent homology barcodes equipped with an $A_\infty$-algebra structure.

The question of stability of this metric is left open by Herscovich, except to note that the 1-ary case coincides with the classical bottleneck distance.

\printbibliography

\end{document}